\magnification=\magstep1
\input amstex
\documentstyle{amsppt}
\TagsOnRight
\NoRunningHeads
\newdimen\notespace \notespace=7mm
\newdimen\maxnote  \maxnote=13mm
\define\marg{\strut\vadjust{\kern-\dp\strutbox
  \vtop to\dp\strutbox{\vss \baselineskip=\dp\strutbox
  \moveleft\notespace\llap{\hbox to \maxnote{\hfil$\to$}}\null}}\ }

\define\cref#1{\marg(#1)}

\redefine\marg{\relax}

\define\dfn#1{{\it #1}}

\define\obs{\strut\vadjust{\kern-\dp\strutbox
  \vtop to\dp\strutbox{\vss \baselineskip=\dp\strutbox
  \moveleft\notespace\llap{\hbox to \maxnote{\hfil!}}\null}}\ }

\define\Dis{\operatorname{Dis}}

\define\id{\operatorname{id}}

\define\Sym{\frak S_n}

\topmatter
\date November 2002\enddate
\title Splitting algebras, Symmetric functions and\\
  Galois Theory
\endtitle 
\author T. Ekedahl and D. Laksov\endauthor
\affil Department of Mathematics,  University of Stockholm, and
Department of Mathematics, KTH, Stockholm\endaffil
\address Department of Mathematics,  S-106 91 Stockholm, Sweden,
resp. S-100 44 Stockholm, Sweden,\endaddress   
\email teke\@matematik.su.se, laksov\@math.kth.se\endemail

\abstract We present a theory for splitting algebras of monic
polynomials over rings, and apply the results to symmetric functions,
and Galois theory. Our main result is that the ring of invariants of a
splitting algebra under the symmetric group almost always is the ring
of coefficients. 
\endabstract

\endtopmatter

\document

\head{Introduction}\endhead

Splitting algebras of monic polynomials with coefficients in a ring
appear in a natural way in many different parts of mathematics, and
the theory of splitting algebras is useful in several context. It
offers, for example, an attractive approach to Galois Theory, and it
forms the natural setting for resultants and discriminants, and for the
study of symmetric polynomials. Moreover, the cohomology rings of flag
spaces are splitting algebras.

 Splitting algebras are constructed in most second courses in algebra,
 although the students and teachers are usually not aware of this (see,
 \cite{B} Definition \S 5 AIV.67, and section IV for an
 exception). For all these reasons it is curious that splitting
 algebras and their properties are unknown to most mathematicians.

 In this article we shall try to compensate for this ignorance by
 presenting two different constructions of splitting algebras of monic
 polynomials with coefficients in a commutative ring with unity, and
 prove the main results about splitting algebras. The most important,
 and surprising, result is that if $f(t)$ is a monic polynomial of
 degree $n$ with coefficients in a ring $A$, then the \dfn{symmetric
 group} $\Sym$ on $n$ letters operates on the splitting algebra $A_f$
 of $f(t)$ in a natural way, such that the ring of invariants under
 this action is $A$, when $2$, or the discriminant of $f(t)$, is
 neither zero nor a zero divisor in $A$. We do not know much about the
 exceptional cases when the ring of invariants is strictly larger than
 $A$. It would be satisfactory to have more information on when this
 happens, and what the non constant invariants look like.

  One of the \dfn{classical} and best known results in algebra is the
 Main Theorem on Symmetric Functions that consists of the following
 three assertions (see e.g. \cite{B} \S 6 Th\'eor\`eme 1 AIV58.):

 \proclaim\nofrills{}Let $t, t_1, t_2, \dots 
  , t_n$  be independent variables over a ring $A$, and let  $s_1,
  s_2, 
  \dots ,    s_n$ be the elementary symmetric functions in these
  variables.
\roster
\item The ring of invariants of $A[t_1, t_2, \dots , t_n]$ under
  permutations of the variables $t_1, t_2, \dots , t_n$ is  $A[s_1,
  s_2, \dots , s_n]$.
\item The elementary symmetric functions $s_1, s_2, \dots , s_n$ are
  algebraically independent over $A$, that is, they are not roots of a
  non-zero polynomial in $n$ variables with coefficients in $A$.
\item The ring $A[t_1, t_2, \dots , t_n]$ is free as a $A[s_1, s_2,
  \dots , s_n]$-module. A basis is given by the monomials
$t_1^{m_1}t_2^{m_2} \cdots   t_n^{m_n}$ with $0\leq m_i< i$
for $i=1, 2, \dots 
  , n$. 
\endroster
\endproclaim

 It is less known that the Main Theorem of Symmetric Functions follows
 from corresponding properties of splitting algebras of \dfn{generic}
 polynomials. The advantage of obtaining the Main Theorem of Symmetric
 Functions from the properties of splitting algebras is that the
 results for splitting algebras are valid in much greater
 generality. As an additional benefit we have that the correspondence
 between the Main Theorem of Symmetric Functions and the properties of
 splitting algebras of generic polynomials make the proofs easy and
 conceptual.

 We also indicate how splitting algebras can be used to obtain an
 attractive approach to Galois Theory. Because of the universal
 character of splitting algebras the \dfn{splitting field} $L$ of a
 polynomial of degree $n$ over a field $K$ is the residue field
 $\varphi:K_f\to L$ of the splitting algebra $K_f$ of $f(t)$ over $K$.
 Via $\varphi$ the Galois Group $\operatorname{G}(L/K)$ of $L$ over
 $K$ is simply the subgroup of the symmetric group $\Sym$ that
 preserve the kernel of $\varphi$. This approach appears to be a
 natural way of considering the Galois Theory of $L$ over $K$,
 essentially equivalent to the original way of considering Galois
 groups as the group of permutations of the roots of a polynomial that
 preserve all the relations between the roots. It is worth noting that
 whereas there are examples in all positive characteristics of
 non-trivial extensions with trivial Galois group, the ring of
 invariants of $K_f$ under the action of $\Sym$ is $K$ when the
 characteristic is different from $2$.

 As a typical illustration of our approach to Galois theory we prove
 that an irreducible polynomial with coefficients in $K$ that has a
 root $f$ in the splitting field $L$ of $f$ has all its roots in
 $L$. We choose an element $\tau$ in $K_f$ such that $\varphi(\tau)
 =l$. By our results on the invariants of $K_f$ under the action of
 $\Sym$ (Theorem \cref{2.4} and Corollary \cref{2.5}) we have that the
 polynomial $\prod_{\sigma\in \Sym} (t-\sigma(\tau))$ has coefficients
 in $K$. Hence the polynomial $h(t) =\prod_{\sigma\in \Sym}
 (t-\varphi\sigma(\tau))$ has coefficients in $K$ and splits
 completely in $L$. The root corresponding to $\sigma
 =\id_{\Sym}$ is $l$. Hence the irreducible polynomial $g(t)$ divides
 $h(t)$ and thus have all its roots in $L$.

 We shall make no effort to systematically exploit the connection
 between splitting algebras and Galois Theory, but content ourselves
 with proving that $L$ is unique up to $K$-isomorphisms, and that the
 field of invariants under the action of $\operatorname{G}(K/L)$ on
 $L$ is equal to $K$ when $f(t)$ is separable over $K$. The proofs are
 of a  different character than those commonly given in
 Galois Theory, and are slightly more technical than those of the rest
 of this article.

\def\sectno{1}

\head{\sectno. Splitting algebras}\endhead

We shall give two different constructions of the splitting algebra of
a monic polynomial with coefficients in a ring. Each of the
constructions reveals a different face of splitting algebras and each
leads to valuable information.

 \subhead{\sectno.1 Notation and conventions}\endsubhead All rings
 will be commutative with unity. We say that an element of a ring $A$
 is \dfn{regular} if it is neither zero, nor a zero-divisor. An
 \dfn{algebra} over a ring $A$ is a ring $B$ with a ring homomorphism
 $\varphi :A\to B$. Often we simply write $\varphi:A\to B$ for an
 $A$-algebra. We say that a ring homomorphism $\chi:B\to C$ between
 two $A$-algebras $\varphi:A\to B$ and $\psi: A\to C$ is an
 $A$-\dfn{algebra homomorphism} when $\psi=\chi\varphi$.

 Let $\chi:B\to C$ be an $A$-algebra homomorphism. For each polynomial 
 $$g(t) =t^n-g_1t^{n-1} +g_2t^{n-2}-\cdots +(-1)^n g_n$$
 in the variable $t$ with coefficients in the $A$-algebra
 $B$  we write
 $$(\chi g)(t) =  t^n -\chi(g_1)t^{n-1} +\chi(g_2) t^{n-2}
 -\cdots +(-1)^n \chi(g_n)$$
 for the corresponding polynomial with coefficients in $C$.

 Let $A$ be a ring and let $t, t_1, t_2, \dots , t_n$ be
 \dfn{independent variables} over $A$. We denote by $A[t, t_1, t_2,
 \dots  , t_n]$ the ring of polynomials in the variables $t, t_1, t_2,
 \dots ,  t_n$ with coefficients in $A$. The \dfn{elementary symmetric
 functions} $s_1, s_2, \dots , s_n$ in the variables $t_1, t_2, \dots ,
 t_n$ are defined by the identity
 $$t^n -s_1t^{n-1} +s_2t^{n-2}-\cdots
 +(-1)^n s_n =(t-t_1)(t-t_2)\cdots (t-t_n)\tag{\sectno.1.1}$$
  in the polynomial ring $A[t, t_1, t_2, \dots , t_n] =A[t_1, t_2,
 \dots  , t_n][t]$.

 \definition{\sectno.2 Definition} Let $f(t) =t^n -f_1t^{n-1} +f_2
 t^{n-2}-\cdots +(-1)^nf_n$ be a \dfn{monic} polynomial with
 coefficients in the ring $A$. The \dfn{splitting algebra} $A_f$ of
 $f(t)$ over $A$ is the \dfn{residue ring} of the polynomial ring
 $A[t_1, t_2, \dots , t_n]$ modulo the ideal generated by the elements
 $s_1-f_1, s_2-f_2, \dots , s_n-f_n$.

 Let 
 $$\varphi_f:A[t_1, t_2, \dots , t_n]\to A_f$$
 be the residue map and write 
 $$\varphi_f(t_i) =\tau_i \quad\text{for}\quad i=1,2, \dots , n.$$ By
 definition we have that $\varphi_f(s_i) =\varphi_f(f_i)$ for $i=1,2,
 \dots , n$. Hence, when we apply $\varphi_f$ to the polynomials on
 both sides of the identity \cref{\sectno.1.1},  we obtain a
 \dfn{complete  splitting}  
 $$t^n -f_1t^{n-1} +f_2t^{n-2} -\cdots +(-1)^n f_n
 = (t-\tau_1)(t-\tau_2)\cdots (t-\tau_n)\tag{\sectno.2.1}$$ 
 of $f(t)$.
 We call the  splitting \cref{\sectno.2.1} the
 \dfn{universal splitting}, and we call $\tau_1, \tau_2, \dots ,
 \tau_n$ the \dfn{universal roots}.
 \enddefinition

 \proclaim{\sectno.3 The universal property}
 For every $A$-algebra $\varphi:A\to B$ and every family of elements
 $\upsilon_1, \upsilon_2, \dots, \upsilon_n$ in $B$ such that we have
 a complete  factorization
 $$ t^n -\varphi(f_1)t^{n-1} +\varphi(f_2)
 t^{n-2}-\cdots +(-1)^n \varphi(f_n)
 =(t-\upsilon_1)(t-\upsilon_2)\cdots (t-\upsilon_n)\tag{\sectno.3.1}$$ 
 of $(\varphi f)(t)$, there 
 is a unique $A$-algebra homomorphism
 $$\psi: A_f =A[\tau_1, \tau_2, \dots , \tau_n] \to B$$
 determined by  $\psi(\tau_i) =\upsilon_i$ for $i=1,2, \dots , n$. 
\endproclaim

 \demo{Proof} Let $\chi:A[t_1, t_2, \dots , t_n] \to B$ be the
 unique $A$-algebra homomorphism determined by $\chi(t_i)
 =\upsilon_i$ for $i=1,2, \dots , n$. Then \cref{\sectno.3.1} can be
 written 
 $$\multline t^n -\chi(f_1) t^{n-1} +\chi(f_2) t^{n-2} -\cdots
 +(-1)^n\chi (f_n)\\
   =(t-\chi(t_1)) (t-\chi(t_2)) \cdots
 (t-\chi(t_n)).\endmultline\tag{\sectno.3.2}$$
 When we apply $\chi$ to the polynomials on both sides of
 \cref{\sectno.1.1} we obtain that 
 $$t^n -\chi(s_1) t^{n-1} +\chi(s_2)t^{n-2} -\cdots +(-1)^n
 \chi(t_n) =(t-\chi(t_1)) (t-\chi(t_2)) \cdots
 (t-\chi(t_n)).\tag{\sectno.3.3}$$ From equations
 \cref{\sectno.3.2} and \cref{\sectno.3.3} it follows that
 $\chi(f_i) =\chi(s_i)$ for $i=1,2, \dots , n$. Hence $\chi$
 factors uniquely through an $A$-algebra homomorphism $\psi:A_f\to B$
 such that $\psi(\tau_i)=\upsilon_i$ for $i=1,2, \dots , n$.  \enddemo

 \example{\sectno.4 Remark} As long as we allow the base ring $A$ to
 be arbitrary some care is needed when we treat uniqueness. It is true
 in general that $t-a$ is regular in $A[t]$ for any $a \in A$ but that
 does not imply that if a monic polynomial has a complete splitting
 then that splitting is unique. As an example, assume that we have $0
 \ne \delta \in A$ with $\delta^2=0$. Then we have that
 $t^2=(t+\delta)(t-\delta)$.  \endexample

 \example{\sectno.5 Uniqueness} It follows from the universal property
 of section \cref{\sectno.3} that $A_f$ is unique up to
 isomorphisms. More precisely, if $\chi:A\to A'$ is another
 $A$-algebra such that we have a complete splitting $(\chi f)(t) =
 (t-\tau'_1) (t-\tau'_2) \cdots (t-\tau'_n)$ with roots $\tau'_1,
 \tau'_2, \dots , \tau'_n$ in $A'$, and such that $A'$ satisfies the
 universal property of \cref{\sectno.3}, then there is a unique
 isomorphism of $A$-algebras $\psi:A_f \to A'$ such that $\psi(\tau_i)
 =\tau'_i$ for $i=1,2, \dots , n$. In fact, the existence of such an
 $A$-algebra homomorphism $\psi$ follows from the universal property
 of $A_f$. On the other hand the universal property of $A'$ guarantees
 the existence of an $A$-algebra homomorphism $\psi':A'\to A_f$ such
 that $\psi'(\tau'_i) =\tau_i$ for $i=1,2, \dots , n$. Finally the
 universality of $A_f$ implies that $\psi'\psi$ is the identity on
 $A$, and the universality of $A'$ that $\psi\psi'$ is the identity on
 $A'$.  \endexample

 \subhead{\sectno.6 Action of the symmetric group}\endsubhead Let
 $\Sym$ be the \dfn{symmetric group} of permutations of the
 integers $\{1,2, \dots , n\}$. For each permutation $\sigma$ we have
 a rearrangement $\tau_{\sigma^{-1}(1)}, \tau_{\sigma^{-1}(2)}, \dots
 , \tau_{\sigma^{-1}(n)}$ of the universal roots of $f(t)$. Hence it
 follows from the universal property \cref{\sectno.3} for $A_f$ that
 there is a corresponding $A$-algebra homomorphism
 $$\varphi_\sigma :A_f \to A_f$$
 determined by $\varphi_\sigma (\tau_i) =\tau_{\sigma^{-1}(i)}$ for
 $i=1,2, 
 \dots , n$. It is clear that for $\rho$ and $\sigma$ in $\Sym$
 we have that
 $$\varphi_{\id_{\Sym}} =\id_{A_f} \quad\text{and that}\quad
 \varphi_{\sigma\rho} =\varphi_\sigma\varphi_\rho.$$
 In other words, the group $\Sym$ \dfn{acts} on the $A$ algebra
 $A_f$. 

\definition{\sectno.7 Alternative construction of the splitting
 algebra} We start the alternative construction with the $A$-algebra
 $A=A_n$ and the polynomial $f(t) =f_n(t) $ and construct, by
 descending induction on $i$, $A$-algebras $A_{i} =A[\upsilon_n,
\upsilon_{n-1}, \dots , \upsilon_{i+1}]$ and polynomials $f_{i}(t)$
in $A_i[t]$ for $i=n-1, n-2, \dots , 0$ as follows:

Assume that we have defined $A_i =A[\upsilon_n, \upsilon_{n-1}, \dots ,
\upsilon_{i+1}]$ and $f_i(t)$ in $A_i[t]$. We define $A_{i-1}$ by
$A_{i-1} =A_i[t]/(f_i(t))$ and we let $\upsilon_i$ be the class of $t$
in $A_{i-1}$. Then $A_{i-1} =A_i[\upsilon_i] =A[\upsilon_n,
\upsilon_{n-1}, \dots , \upsilon_i]$, and we define $f_{i-1}(t)$ by
$f_{i-1}(t) =f_i(t)/(t-\upsilon_i)$.

Note that $A_0 =A_1$ since $f_1(t) =t-\upsilon_1$.

The $A$-algebra $A_0$ is canonically isomorphic to $A_f$. In fact, the
$A$-algebra $A_0$ has the universal property \cref{\sectno.3}.
To see this
we note that $f(t)$ splits as $f(t) =\prod_{i=1}^n (t-\upsilon_i)$ in
$A_0[t]$, and for every $A$-algebra $\varphi:A\to B$ such that
$(\varphi f)(t) =\prod_{i=1}^n (t-\sigma_i)$ with $\sigma_1, \sigma_2,
\dots , \sigma_n$ in $B$ we can, by descending induction on $i$,
construct an $A$-algebra homomorphism $A_i \to B$ for $i=n, n-1, \dots
, 0$ which maps $\upsilon_j$ to $\sigma_j$ for $j=n, n-1, \dots , i$.
\enddefinition

\definition{\sectno.8 Remark}A useful consequence of the alternative 
construction is that the natural map $A[\tau_n]_{f_{n-1}} \to A_f$ is
an $A$-algebra isomorphism.
\enddefinition

\definition{\sectno.9 Remark} The formation of splitting algebras is
compatible with \dfn{extension} of scalars. That is, when $\varphi:
A\to B$ is a ring homomorphism the canonical map $A_f\to B_{\varphi
f}$ resulting from the universal splitting of $\varphi f$ in
$B_{\varphi f}$ gives, by extension of scalars, an isomorphism of
$B$-algebras
 $$B\otimes _A A_f \to B_{\varphi f}.$$
 \enddefinition

\proclaim{\sectno.10 Proposition}The splitting algebra $A_f =A[\tau_1,
 \tau_2, \dots , \tau_n]$ is a free module over the ring $A$ with a
basis consisting of the elements $\tau_2^{m_2} \tau_3^{m_3} \cdots
\tau_n^{m_n}$ with $0\leq m_i <i$ for $i=2,3, \dots , n$.
\endproclaim

 \demo{Proof}The polynomial $f_i(t)$ of the alternative construction
 \cref{\sectno.7} is monic of degree $i$, and has coefficients in
 $A$. Hence the $A_i$-algebra $A_{i-1} =A_i[t]/ (f_i(t))$ is a free
 $A_i$-module of rank $n$ with basis $1, \upsilon_i, \dots ,
 \upsilon_i^{n-1}$ for $i=1,2, \dots , n$. The proposition
 consequently follows from the inductive definition of $A_0=A_1$.
 \enddemo

 \def\sectno{2}

\head{2. Splitting algebras of generic polynomials}\endhead 

We shall prove the main properties of splitting algebras of
\dfn{generic} polynomials, that is, polynomials whose coefficients are
independent variables over a base ring. These properties are collected
in the Main Theorem on Splitting Algebras of Generic Polynomials. The
properties are analogous to the properties of the Main Theorem of
Symmetric Functions. In fact we observe that the two \dfn{Main
Theorems}  are equivalent.

 \definition{\sectno.1 Definition}Let $f(t) =t^n -f_1t^{n-1}
 +f_2t^{n-2} -\cdots +(-1)^n f_n$ be a polynomial in the variable $t$
 with coefficients in  the ring $A$, and let $A_f =A[\tau_1, \tau_2,
 \dots , \tau_n]$ be the splitting algebra of $f(t)$ over $A$ with
 universal roots $\tau_1, \tau_2, \dots , \tau_n$. We write
 $$\Dis(f) =(-1)^{n(n-1)/2} \prod_{i\neq j} (\tau_i-\tau_j)
 =\prod_{i>j} (\tau_i-\tau_j)^2$$
 and call $\Dis(f)$ the \dfn{discriminant} of $f(t)$.
 \enddefinition

\proclaim{\sectno.2 Lemma}Let $B$ be an $A$-algebra and let $f(t) =t^n
 -f_1 t^{n-1} +f_2 t^{n-2}-\cdots +(-1)^n f_n$ be a polynomial with
 coefficients 
in $B$. Write $B_f =B[\tau_1, \tau_2, \dots , \tau_n]$ where $\tau_1,
\tau_2, 
\dots , \tau_n$ are the roots of $f(t)$ in $B_f$. Assume that the
discriminant $\Dis(f)$ is regular in
$B_f$. Then $B$ is the ring of invariants of $B_f$ under the action of
$\Sym$.
\endproclaim

\demo{Proof} Let $\upsilon$ be an element in $B_f$ which is invariant
under the action of $\Sym$. It follows from Proposition \cref{1.10}
that we have a unique expression
$$\upsilon=\sum_{0\leq m_i<i, s\leq i\leq n} g(m_s, m_{s+1}, \dots ,
m_n)\tau_s^{m_s} \tau_{s+1}^{m_{s+1}} \cdots \tau_n^{m_n},$$ with
$g(m_s, m_{s+1}, \dots , m_n)$ in $B$ and where $s\geq 2$. The lemma
asserts that all the coefficients $g(m_s, m_{s+1}, \dots , m_n)$ are
zero except $g(0,0, \dots , 0)$. Assume to the contrary that at least
one of the coefficients $g(m_s, m_{s+1}, \dots , m_n)$ with $m_s>0$ is
non-zero.  Since $\upsilon$ is invariant under $\Sym$ it is invariant
under the transposition $\sigma$ which interchanges $\tau_{s-1}$ and
$\tau_s$ and keeps the remaining roots fixed. We obtain that
 $$(\upsilon-\sigma(\upsilon)) =\sum_I g(m_s,
 m_{s+1}, \dots , m_n)(\tau_s^{m_s} 
 -\tau_{s-1}^{m_s})( \tau_{s+1}^{m_{s+1}} \cdots \tau_n^{m_n})
 =0,\tag{\sectno.2.1}$$  
 where the sum is over the set $I= \{0\leq m_i<i, s\leq i\leq n,
   m_s>0\}$. 
 By assumption the element $\prod_{i\neq j}(\tau_i-\tau_j)$, and thus
 each element  $(\tau_i-\tau_j)$, is regular in
 $B_f$. Hence  we can divide the term $(\tau_s^{m_s}
 -\tau_{s-1}^{m_s})$  of (\sectno.2.1) by 
 $(\tau_s-\tau_{s-1})$ to obtain the equation
 $$\sum_I g(m_s, m_{s+1}, \dots, m_n)(\sum_{j=1}^{m_s-1} \tau_{s-1}^j
 \tau_s^{m_s-1-j})(\tau_{s+1}^{m_{s+1}} \cdots
 \tau_n^{m_n})=0.\tag{\sectno.2.2}$$
 All the monomials
 $\tau_{s-1}^j\tau_s^{m_s-1-j} \tau_{s+1}^{m_{s+1}} \cdots
 \tau_n^{m_n}$ 
 on the left hand side of  equation (\sectno.2.2) are different and
 none of them contains a power of $\tau_1$ since, if $s=2$ and $m_s>0$,
 we have that $m_s-1 =m_2-1=0$. Hence \cref{\sectno.2.2} is a
 non-trivial  relation 
 between the monomials $\tau_2^{m_2} \tau_3^{m_3} \cdots \tau_n^{m_n}$
 with $0\leq m_i <i$. This contradicts the assertion of Proposition
 \cref{1.10}, and thus contradicts the assumption that some
 coefficient 
 $g(m_s,  m_{s+1}, \dots , m_n)$ with $m_s>0$ is non-zero.  \enddemo

\proclaim{\sectno.3 Lemma}Let $t_1, t_2, \dots , t_n$ be algebraically
independent elements over $A$, and let $s_1, s_2, \dots , s_n$ be the
elementary symmetric functions in these variables. Moreover, let $f(t)
=t^n -f_1t^{n-1}+ f_2t^{n-2}-\cdots+(-1)^nf_n$ be a polynomial whose
coefficients 
are algebraically independent over $A$, and let $\tau_1, \tau_2, \dots
, \tau_n$ be the roots of $f(t)$ in the splitting algebra $A[f_1, f_2,
\dots 
, f_n]_f$ for $f(t)$ over $A[f_1, f_2, \dots , f_n]$. Then there is an
isomorphism of $A$-algebras
 $$A[\tau_1, \tau_2, \dots , \tau_n] \to A[t_1, t_2, \dots , t_n]$$
 mapping $\tau_i$ to $t_i$, and thus $f_i$ to $s_i$, for $i=1, 2,
 \dots ,  n$. 
\endproclaim

\demo{Proof}Since $f_1, f_2, \dots , f_n$ are algebraically
independent over $A$ by assumption, there is an $A$-algebra
homomorphism $\varphi: A[f_1, f_2, \dots , f_n] \to A[s_1, s_2, \dots
, s_n]$ 
mapping the variable $f_i$ to $s_i$ for $i=1,2, \dots , n$. We obtain
that 
$(\varphi f)(t)= t^n -s_1t^{n-1}+s_2t^{n-2}- \cdots +(-1)^n s_n=
\prod_{i=1}^n (t-t_i)$.  By the 
universal property of  splitting algebras there is a unique
$A$-algebra homomorphism $\psi: A[\tau_1, \tau_2, \dots , \tau_n]=
A[f_1, 
f_2,
\dots , f_n]_f \to A[t_1, t_2, \dots , t_n]$ such that $\psi(\tau_i)
=t_i$ for $i=1, 2, \dots , n$, and restricting to $\varphi$ on $A[f_1,
f_2, \dots, f_n]$.  

We shall construct an inverse to this homomorphism. By assumption the
elements $t_1, t_2, \dots , t_n$ are algebraically independent over
$A$. Consequently there is an $A$-algebra homomorphism $A[t_1, t_2,
\dots , t_n] \to A[\tau_1, \tau_2, \dots , \tau_n]$ mapping $t_i$ to
$\tau_i$ for all $i$. It is clear that this homomorphism is the
inverse to $\psi$.
\enddemo

\proclaim{\sectno.4 The Main Theorem on Splitting Algebras of Generic
Polynomials}Let
$f(t) =t^n -f_1 t^{n-1} + f_2t^{n-2}-\cdots +(-1)^nf_n$ be a
polynomial whose 
coefficients are algebraically independent over the ring $A$. The
following three assertions hold:
\roster
\item The ring of invariants of the  splitting algebra
  $A[f_1, f_2, \dots , f_n]_f$ of $f(t)$ over $A[f_1, f_2, \dots ,
  f_n]$ 
  under the action of\/ $\Sym$ is $A[f_1, f_2, \dots , f_n]$.
\item The universal roots $\tau_1, \tau_2, \dots , \tau_n$ of the
  polynomial $f(t)$  in $A[f_1, 
  f_2, \dots 
  , f_n]_f$  are algebraically independent over
  $A$.
\item The ring $A[f_1, f_2, \dots , f_n]_f$ is free as an $A[f_1, f_2,
  \dots , f_n]$-module. A basis is given by the monomials
  $\tau_2^{m_2}\tau_3^{m_3} \cdots \tau_n^{m_n}$ with $0\leq m_i<i$ for
  $i=2,3, \dots , n$.
\endroster
\endproclaim

\demo{Proof}It follows from Lemma \cref{\sectno.3} that $\tau_1,
  \tau_2, \dots , \tau_n$ are algebraically independent over $A$. That
is, the second property of the theorem holds. In particular
$\prod_{i\neq j} (\tau_i-\tau_j)$ is neither zero nor a zero divisor
in $A_f =A[\tau_1,
\tau_2, \dots , \tau_n]$. Hence the first property of the theorem
follows from Lemma \cref{\sectno.2}. The third property follows from
 Proposition \cref{1.10}.
 \enddemo

 \proclaim{\sectno.5 Corollary} Let $f(t) =t^n -f_1t^{n-1} +f_2t^{n-2}
 -\cdots +(-1)^n f_n$ be a polynomial with coefficients in $A$ and let
 $\tau_1, \tau_2, \dots, \tau_n$ be the universal roots of $f$ in
 $A_f$. Moreover  let
 $\varphi:A_f \to B$ be a homomorphism of $A$-algebras. For all
 polynomials $h(t_1, t_2, \dots , t_n)$ in the variables $t_1, t_2,
 \dots , t_n$ with coefficients in $A$ that are invariant under the
 action of $\frak S_n$ on the variables, we have that
 $h(\varphi(\tau_1), \varphi(\tau_2,), \dots , \varphi(\tau_n))$ is in
 $A$.

 In particular the discriminant $\Dis(f)$ lies in $A$, and if \/
 $\Dis(f)$ is regular in  $A$ we have that $A$ is the ring of
 invariants of   $A_f$ under the action of $\Sym$.
 \endproclaim

  \demo{Proof} Let $g(t) = t^n -g_1t^{n-1}+f_2t^{n-2}-\cdots+(-1)^n
  g_n$ be a polynomial whose coefficients are algebraically
  independent over $A$, and let $A_g =A[\upsilon_1, \upsilon_2, \dots
  , \upsilon_n]$ be the generic splitting algebra, where $\upsilon_1,
  \upsilon_2, \dots , \upsilon_n$ are the universal roots of
  $g(t)$. By the universal property of splitting algebras the
  homomorphism $A[g_1, g_2, \dots , g_n]\to A$ that maps $g_i$ to
  $f_i$ for $i=1,2, \dots , n$ can be extended to an $A$-algebra
  homomorphism $\psi:A[g_1, g_2, \dots , g_n] \to A_f$ determined by
  $\psi(\upsilon_i) =\tau_i$ for $i=1,2, \dots , n$.

 It follows from assertion (1) Theorem \cref{2.4} that $h(\upsilon_1,
  \upsilon_2, \dots , \upsilon_n)$ is in $A$, and consequently that
  $$\multline
 h(\varphi(\tau_1) , \varphi(\tau_2), \dots ,\varphi(\tau_n)) 
  =\varphi(h(\tau_1, \tau_2, \dots , \tau_n))\\
  =\varphi(h(\psi(\upsilon_1), \psi(\upsilon_2), \dots ,
  \psi(\upsilon_n)) =\varphi\psi(h(\upsilon_1, \upsilon_2,
 \upsilon_n)\endmultline$$
  is in $A$. Hence we have proved the first part of the Corollary.

 When $h(t_1, t_2, \dots, t_n) =(-1)^{n(n-1)/2} \prod_{i\neq j} 
  (t_i-t_j)$, and $B=A_f$ with $\varphi=\id_{B}$, we obtain that
  $h(\varphi(\tau_1), \varphi(\tau_2), \dots , \varphi(\tau_n))
  =h(\tau_1, \tau_2, \dots , \tau_n) =\Dis(f)$. Hence the second part
  of the Corollary is a consequence of the first part.

 It follows from Proposition \cref{1.10} that if $\Dis(f)$ is regular
in $A$ then it is regular in $A_f$. Hence $A$ is the ring of
invariants of $A_f$ under $\Sym$ when $\Dis(f)$ is regular
in $A$. Hence we have proved the last part of the Corollary.
\enddemo

\proclaim{\sectno.6 Theorem} The Main Theorem on Symmetric Functions
is equivalent to the Main Theorem on Splitting Algebras of Generic
Polynomials.  \endproclaim
  
  \demo{Proof}The assertion of the theorem follows immediately from
  Lemma \cref{\sectno.3} since the algebraic independence of the
  elements  $t_1, t_2,
  \dots , t_n$ implies that the elements $\tau_1, \tau_2, \dots ,
  \tau_n$ 
  are algebraically independent, and  the algebraic independence
  of the elements $f_1, f_2, \dots , f_n$ implies that the elements
  $s_1, s_2,   \dots , s_n$ are 
  algebraically independent.  \enddemo

  \definition{\sectno.7 Remark}It is curious that it is the
  algebraic independence of the elements $f_1, f_2, \dots , f_n$ that
  implies the algebraic independence of $s_1, s_2, \dots , s_n$,
  whereas it is the algebraic independence of the elements $t_1, t_2,
  \dots , t_n$ that implies the algebraic independence of $\tau_1,
  \tau_2, \dots , \tau_n$.  This fact, accounts for the difference in
  the second properties of the Main Theorem on Symmetric Functions and
  the Main Theorem on Splitting Algebras of Generic Polynomials.
 \enddefinition

\def\sectno{3}  

\head{\sectno. The ring of invariants of splitting algebras}\endhead

 We shall now give a different proof of Lemma\cref{2.2}. As will be
 seen it also gives the conclusion of the Lemma under a different
 condition, namely that $2$ is regular in $A$. It should be noted that
 the proof is completely independent of the previous section and thus
 gives an alternative proof of Lemma\cref{2.2} and of the first
 assertion of Theorem\cref{2.4}.

\proclaim{\sectno.1 Theorem} Let $f(t) =t^n -f_1t^{n-1}+f_2t^{n-2}-
 +\cdots +(-1)^n 
f_n$ be a polynomial with coefficients in the ring $A$. Assume that
either of the following two conditions hold:
\roster
\item The element $2$ is regular in $A$.
\item The element $\prod_{i\neq j} (\tau_i -\tau_j)$ is regular in
$A_f$. 
\endroster
Then $A$ is the ring of invariants of
$A_f$ under the action of $\Sym$.
\endproclaim

\demo{Proof} When the degree $n$ of $f(t)$ is  $1$ the assertions of
the theorem obviously hold. We shall prove the case $n=2$. It follows
from Proposition \cref{1.10} that an element which is invariant under
$\frak S_2$ is of the form $g+h\tau_1 =g+h\tau_2$ with $g$ and $h$ in
$A$.  Thus we must have that $h(\tau_1 -\tau_2) =h(f_1-2\tau_2)=0$. In
particular it follows from Proposition \cref{1.10} that $2h=0$. When
$(\tau_1-\tau_2)(\tau_2-\tau_1)$, or $2$, are regular, we obtain that
$h=0$. Hence the theorem holds for $n=2$.

We now prove the theorem by induction on the degree of $f(t)$,
starting with the already proved cases of degree $1$ and $2$. Assume
that $n\geq 3$ and that the theorem holds for $n-1$. Write $A_f
=A[\tau_1, \tau_2, \dots , \tau_n]$ where $\tau_1, \tau_2, \dots ,
\tau_n$ are the roots of $f(t)$. It follows from Remark
\cref{1.8} that $A_f$ is the splitting algebra of $f_{n-1}(t)
=f(t)/(t-\tau_n)$ over $A[\tau_n]$. If $2$ is regular in $A$ it
follows that it is regular in $A[\tau_n]$. Moreover, if $\prod_{i\neq
j} (\tau_i-\tau_j)$ is regular in $A_f$, then $\prod_{i\neq j}
(\tau_i-\tau_j)$ is regular in $A[\tau_n]_{f_{n-1}} =A_f$. By the
induction hypothesis applied to the splitting algebra
$A[\tau_n]_{f_{n-1}}$ of $f_{n-1}$ over $A[\tau_n]$ every invariant in
$A_f$ under $\Sym$ will be in $A[\tau_n]$ and can therefore be written
as $g(\tau_n)$ for some polynomial $g(t)$ of degree strictly less that
$n$ and with coefficients in $A$.

The element $g(\tau_n)$ is invariant under the transposition
exchanging $\tau_{n-1}$ and $\tau_n$ and leaving the remaining roots
fixed. Hence we have that $g(\tau_n) =g(\tau_{n-1})$ in $A[\tau_n,
\tau_{n-1}] =A[\tau_n][t]/(f_{n-1}(t))$. Consequently we have that
$g(t) -g(\tau_n) =h(t) f_{n-1}(t)$ in $A[\tau_n][t]$ for some
polynomial $h(t)$. Since $g(t)$ is of degree strictly less than $n$ we
see that $h(t)$ is of degree $0$ in $t$. Hence we have that $h
=k(\tau_n)$ where $k(s)$ is a polynomial in the variable $s$ of degree
strictly less than $n$ with coefficients in $A$. Hence we have that
$g(t) -g(\tau_n) =k(\tau_n) f_{n-1}(t)$ for some polynomial $k(t)$ in
$A[\tau_n][t]$. Multiplication by $(t-\tau_n)$ gives that
 $$(t-\tau_n)(g(t) -g(\tau_n)) =k(\tau_n)f(t)$$
 in $A[\tau_n][t] =A[t][s]/(f(s))$. Consequently we have that
 $$(s-t) (g(s) -g(t)) =k(s)f(t) +l(s,t)f(s)\tag{\sectno.1.1}$$ for
 some polynomial $l(s.t)$ in $A[s,t]$. The left hand side of
 \cref{\sectno.1.1} has degree at most $n$ in $s$ and $k(s)$ is of
 degree strictly less than $n$ in $s$. Hence $l(s,t)$ must be equal to
 the coefficient $g$ of $s^n$ in $s g(s)$. With $s=t$ in
 \cref{\sectno.1.1} we see that $k(s) =-l(s,s) =g$.  We have proved
 that
 $$(s-t)(g(s) -g(t)) =g(f(s) -f(t)).\tag{\sectno.1.2}$$
 Comparing the homogeneous component of degree $n$ on both sides of
 (\sectno.1.2)  we  get
 $$(s-t)(gs^{n-1} -gt^{n-1}) =g(s^n -t^n)\tag{\sectno.1.3}$$
 in $A[s,t]$.
 The equation (\sectno.1.3) can not be satisfied when $n\geq 3$ unless
 $g=0$ and thus $(s-t)(g(s) -g(t))=0$. Consequently $g(s) =g(t)$ so
 that $g(t)$ is in $A$, and so is $g(\tau_n)$, as we wanted to prove.
 \enddemo

\example{\sectno.2 Example}When condition (1) of Theorem
 \cref{\sectno.1} does not hold, the conclusion of the theorem may not
 be true. For example, the elements of the form $g+h\tau_1$ in the
 splitting algebra $A[\tau_1, \tau_2]$ of the polynomial $f(t) =t^2
 -f_1t +f_2$ are invariant under $\frak S_2$ when $h$ is in the
 annihilator of $f_1$ and $2h=0$.  This is because $g+h\tau_1
 =g+h\tau_2$ exactly when $h(\tau_1-\tau_2) =h(f_1 -2\tau_2)=0$. In
 particular, when $2=0$ in $A$ and $A$ has zero divisors the ring of
 invariants of $A_f$ under the action of $\frak S_2$ is, for a
 suitable polynomial $f(t)$, strictly bigger than $A$.  \endexample

\example{\sectno.3 Remark}We saw in Corollary \cref{2.5} that the
discriminant of $f(t)$ is in $A$, and that we can replace the second
condition in Theorem \cref{\sectno.1} by the condition that the
discriminant does not divide zero in $A$. We have chosen to present
Theorem \cref{\sectno.1} in the apparently less general form to
emphasize that the proof is independent of the results of Section 2.
\endexample

 \def\sectno{4}

\head{\sectno. Some uses of splitting algebras to Galois
theory}\endhead 

 We shall illustrate some of the uses of splitting algebras to Galois
 theory. The first two results hold over an arbitrary base ring
 $A$. We say that two polynomials $f(t)$ and $g(t)$ in $A[t]$ are
 \dfn{mutually prime} when the ideal they generate is equal to
 $A[t]$. 

\proclaim{\sectno.1 Lemma}Let $f(t)=\prod_{i=1}^r g_i(t)$ be a monic
polynomial of degree $n$ which is the product of monic polynomials
$g_i(t)$ of degree $n_i$ that are mutually prime over the field
$A$. Denote by $A_{g_i} =A[\upsilon_{1i},
\upsilon_{2i}, \dots ,
\upsilon_{n_ii}]$ the splitting algebra of $g_i(t)$ over $A$ where the
elements $\upsilon_{1i}, \upsilon_{2i}, \dots , \upsilon_{n_ii}
=\upsilon_i$ are the universal roots of $g_i(t)$, and let $A_f
=A[\tau_1, \tau_2, 
\dots , \tau_n]$, where $\tau_1, \tau_2, \dots , \tau_n$ are the
universal roots of $f(t)$. Let $h_i(t) =g_1(t)\cdots
g_i(t)/(t-\upsilon_i)\cdots g_r(t)$ in 
$A[\upsilon_i][t]$. Then we have a canonical isomorphism
 $$A_f \to \prod_{i=1}^r A[\upsilon_i]_{h_i}.$$
\endproclaim

\demo{Proof}It follows from the Chinese Remainder Theorem that we have
a canonical isomorphism $A[t]/(f) \to \prod_{i=1}^r A[t]/(g_i)$, that
is 
$$A[\tau_n] \to \prod_{i=1}^r A[\upsilon_i].$$ The image of
$f(t)/(t-\tau_n)$ by the projection $A[\tau_n][t] \to A[\upsilon_i][t]$
is $h_i$. It is clear that since $g_1(t), g_2(t), \dots , g_r(t)$ are
mutually 
prime over $A$, the polynomials $g_1(t), \dots , h_i(t), \dots ,
g_r(t)$ are 
mutually prime over $A[\upsilon_i]$. From the universal property of
splitting algebras it consequently follows that
 $$A_f = A[\tau_{n}]_{f/(t-\tau_n)} =\prod_{i=1}^r
 A[\upsilon_i]_{h_i},$$ which is the isomorphism of the lemma.
 \enddemo

 \example{\sectno.2 Notation}Let $n_1, n_2, \dots , n_r$ be positive
integers and let $n=n_1 +n_2+\cdots +n_r$. An $(n_1, n_2, \dots ,
n_r)$-{\it shuffle} of the integers $1,2, \dots , n$ is a permutation
$\sigma$ in $\Sym$ such that when $n_0=0$ and $p_i =n_1+n_2+\cdots
+n_j$ for $j=1,2, \dots , r$, we have that
 $$\sigma (p_{j-1}+1)< \sigma(p_{j-1}+2)< \cdots < \sigma(p_j),$$
 for $j=1,2, \dots , r$. The set $\frak S(n_1 \vert n_2\vert \cdots
 \vert n_r)$  of $(n_1, n_2, \dots , n_r)$ shuffles of $1,
 2, \dots , n$ represents the left cosets in $\Sym$ of the subgroup
 $\frak S(n_1, n_2, \dots , n_r)$ that permutes the elements in the
 sets $(p_{j-1}+1, p_{j-1}+2, \dots , p_j)$ for $j=1, 2, \dots , r$.

Let $f(t)=\prod_{i=1}^r g_i(t)$ be a product of polynomials $g_i(t)$
of degree 
$n_i$. For each shuffle $\sigma$ in $\frak S(n_1\vert n_2\vert\cdots
\vert n_r)$ we obtain a map 
 $$\varphi_\sigma \colon A_f \to A_{g_1}\otimes_A A_{g_2}\otimes _A
 \cdots \otimes _A A_{g_r}$$ which maps
 $\tau_{\sigma^{-1}(p_{j-1}+1)}, \tau_{\sigma^{-1}(p_{j-1}+2)}, \dots
 , \tau_{\sigma^{-1} (p_j)}$ to $\rho_{1j} =1\otimes_A \cdots
 \otimes_A \upsilon_{1j}\otimes_A \cdots \otimes_A1, \rho_{2j} =
 1\otimes_i \cdots \otimes_A \upsilon_{2j}\otimes_A \cdots \otimes_A
 1, \dots, \rho_{n_jj} =1\otimes_A \cdots \otimes_A \upsilon_{n_jj}
 \otimes_A \cdots \otimes_A 1$ for $j=1, 2, \dots , r$ where
 $\upsilon_{ij}$ is in the $j$'th factor.  The map $\varphi_\sigma$
 maps $\tau_1, \tau_2, \dots , \tau_n$ to the roots $\rho_{11},
 \rho_{21}, \dots , \rho_{n_11}, \rho_{12}, \rho_{22} , \dots ,
 \rho_{n_22} , \dots , \rho_{1r}, \rho_{2r}, \dots , \rho_{n_rr}$
 shuffled by $\sigma$. That is $\varphi_\sigma(\tau_i)
 =\rho_{\sigma(A)j}$ where $j$ is determined by $p_{j-1} <i \leq p_j$.
  \endexample

 \proclaim{\sectno.3 Theorem}Let $f(t)=\prod_{i=1}^r g_i(t)$ be a
 splitting over $A$ of the monic polynomial $f(t)$ in mutually prime
 monic polynomials $g_i(t)$. Then  the homomorphism
 $$\prod_{\sigma \in \frak S(n_1\vert n_2\vert \cdots \vert n_r)}
 \varphi_\sigma \colon A_f
 \to \prod_{\sigma \in \frak S(n_1\vert n_2\vert \cdots \vert n_r)}
 A_{g_1}\otimes_A A_{g_2}\otimes_A \cdots \otimes_A A_{g_r}$$
 is an isomorphism.
\endproclaim

\demo{Proof}We prove the theorem by induction on the degree $n$ of
$f(t)$. It is clear that the assertion of the theorem holds when
$n=1$. Assume that the theorem holds for polynomials of degree
$n-1$. It follows from Lemma \cref{\sectno.1} that we have an
isomorphism
 $$A_f \to \prod_{i=1}^r A[\upsilon_i]_{h_i}$$
 with $h_i(t) =g_1(t)\cdots g_i(t) /(t-\upsilon_i)\cdots g_r(t)$. From
 the 
 induction assumption it follows that there is an isomorphism
 $$A[\upsilon_i] \to \prod_{\sigma \in \frak S_i}
 A[\upsilon_i]_{g_i}\otimes_{A[\upsilon_i]} \cdots\otimes
 _{A[\upsilon_i]}
 A[\upsilon_i]_{g_i/(t-\upsilon_i)}\otimes_{A[\upsilon_i]} \cdots
 \otimes_{A[\upsilon_i]} A[\upsilon_i]_{g_r}$$ where the product is
 over all the $(n_1, \dots , n_i-1, \dots, n_r)$-shuffles $\frak S_i$
 of the integers $(1, 2, \dots , \dots, n_i-1 , n_i+1, \dots ,
 n_r)$. We have that $A[\upsilon_i]_g =A[\upsilon_i]\otimes_A A_g$ for
 all polynomials $g$ in $A[t]$, and that
 $A[\upsilon_i]_{g_i/(t-\upsilon_i)} =A_{g_i}$. Consequently it
 follows from the induction hypothesis that $A[\upsilon_i]_{h_i}
 =\prod_{\sigma\in \frak S_i} A_{g_1} \otimes_A A_{g_2} \otimes_A\cdots
 \otimes_A A_{g_r}$. We thus have an isomorphism
 $$A_f \to \prod_{i=1}^r A[\upsilon_i]_{h_i} \to \prod_{i=1}^r
 \prod_{\sigma\in 
\frak S_i} A_{g_1}\otimes_A A_{g_2}\otimes_A \cdots \otimes_A
A_{g_r}.\tag{\sectno.3.1}$$
The right hand side of \cref{\sectno.3.1} is canonically isomorphic to
$\prod_{\sigma\in 
  \frak S(n_1\vert n_2\vert\cdots \vert n_r)} A_{g_1}\otimes_A A_{g_2}
\otimes_A \cdots \otimes_A A_{g_r}$ by the map that maps
  $\tau_n$ to $\rho_{n_ii}$ and coincides with $\sigma$ on $\tau_1,
  \tau_2, \dots , \tau_{n-1}$. Consequently the isomorphism
  \cref{\sectno.3.1} is the map described   in the theorem.
\enddemo

\example{\sectno.4 Remark}A consequence of Theorem \cref{\sectno.3} is
  that $A_f$ as a $\Sym$-representation is induced from the $\frak
  S_{\mu}$-represen\-tation $A_{g_1}\bigotimes\cdots\bigotimes
  A_{g_r}$, where $\mu$ is the partition $\{n_1,\dots,n_r\}$.
\endexample

\proclaim{\sectno.5 Proposition}Let $f(t)$ be a monic polynomial with
coefficients in a field $K$. The group $\Sym$ operates
transitively on the maximal ideals of $K_f$.
\endproclaim

\demo{Proof}Let $L$ be an algebraic extension of $K$ where the
polynomial $f(t)$ splits. The maximal ideals in $K_f$ are the
intersection of the maximal ideals in $L_f =L\otimes_K K_f$ by the
integral extension $K_f \to L\otimes_K K_f$. We can therefore assume
that $f(t)$ splits in linear factors in $K$.

Write $f(t)= \prod _{i=1}^r (t-f_i)^{n_i}$, with $f_1, f_2, \dots ,
f_r$ different elements in $K$. It follows from Theorem
\cref{\sectno.3} that we have $K_f =\prod_{\sigma\in \frak S(n_1 \vert
n_2\vert \cdots
\vert n_r)} K_{(t-f_1)^{n_1}} \otimes_K K_{(t-f_2)^{n_2}}\otimes_K
K_{(t-f_r)^{n_r}}$. Consequently every maximal ideal in $K_f$ comes
from a maximal ideal in one of the factors $K_{(t-f_1)^{n_1}}
\otimes_K K_{(t-f_2)^{n_2}}\otimes_K \cdots \otimes_K
K_{(t-f_r)^{n_r}}$. For each positive integer $m$ and each element $g$
of $K$ we have a $K$-algebra homomorphism $K_{(t-g)^m} \to
K_{(t-g)}=K$ which maps the roots $\upsilon_1, \upsilon_2, \dots ,
\upsilon_m$ of $(t-g)^m$ in $K_{(t-g)^n}$ to $g$.  The kernel of
$K_{(t-g)^m} \to 
K_{(t-g)}$ is therefore the maximal ideal generated by $\upsilon_1-g,
\upsilon_2-g, \dots , \upsilon_r-g$, and it is nilpotent because
$(\upsilon_i-g)^m=0$ for all $i$.  Hence $K_{(t-f_i)^{n_i}}$ has a
single maximal ideal $\frak I_i$ whose residue field is $K$. It
follows that 
$K_{(t-f_1)^{n_1}} \otimes_K K_{(t-f_2)^{n_2}}\otimes_K
K_{(t-f_r)^{n_r}}$ has a single maximal ideal $\sum_{i=1}^r
K_{(t-f_1)^{n_1}} \otimes_K \cdots \otimes_K \frak I_i\otimes_K \cdots
\otimes_K K_{(t-f_r)^{n_r}}$ whose residue field is $K$. Since $\Sym$
operates transitively on the factors of $\prod_{\sigma\in \frak
S(n_1\vert n_2\vert \cdots \vert n_r)} K_{(t-f_1)^{n_1}} \otimes_K
K_{(t-f_2)^{n_2}}\otimes_K K_{(t-f_r)^{n_r}}$ we have that it operates
transitively on the maximal ideals.  \enddemo

\proclaim{\sectno.6 Corollary}A splitting field of a polynomial with
coefficients in a field $K$ is determined up to $K$-isomorphisms.
\endproclaim

\demo{Proof}Let $L$ be a splitting field of the polynomial $f(t)$ with
coefficients in $K$. We have a natural surjection $K_f =K[\tau_1,
\tau_2, \dots , \tau_n] \to L$ mapping the roots $\tau_1, \tau_2,
\dots , 
\tau_n$ of $f(t)$ in $K_f$ to the roots of $f$ in $L$, in some
order. The 
kernel is a maximal ideal. Hence the corollary follows from the
proposition.
\enddemo

\proclaim{\sectno.7 Proposition}Let $f(t)$ be a polynomial which is
separable over a field $K$, that is $f(t)$ and its formal derivative
$f'(t)$ are mutually prime over $K$. Then we have that \roster \item
The 
field $K$ is the ring of invariants of $K_f$ under the action of
$\Sym$.
\item The splitting algebra $K_f$ can be written as a direct product
$K_f =\prod_{i=1}^s K_i$ of fields $K_i$ that are separable over $K$.
\endroster
 In particular the action of $\Sym$ on $K_f$ permutes the factors in
 the product $\prod_{i=1}^s K_i$ transitively, and all the $K_i$ are
 $K$-isomorphic to the same field $L$. Moreover, if $G$ is the
 stabilizer of one of the factors, then $G$ operates on $L$, and $K$
 is the ring of invariants under this action.
\endproclaim

\demo{Proof}Since $f(t) =\prod_{i=1}^n (t-\tau_i)$ over $K_v$ we have
 that $f'(\tau_j) =\prod_{i\neq j} (\tau_j-\tau_i)$. Hence $\Dis(f)
 =(-1)^{n(n-1)/2} \prod_{j=1}^n f'(\tau_j)$. When $f(t)$ is separable
 over $K$  we therefore have that $\Dis(f)$ is invertible in
$K_f$. It follows from Lemma \cref{2.2} or
Theorem \cref{3.1} that the ring of invariants of $K_f$ under the
action of $\Sym$ is $K$.

We prove the second part of the proposition by induction on the degree
$n$ of $f(t)$. The proposition holds for polynomials of degree $1$.
Assume that it holds for polynomials of degree strictly smaller than
$n$.

Write $f(t)$ as a product $f(t) = \prod_{i=1}^r g_i(t)$ of irreducible
polynomials $g_i(t)$ in $K[t]$ that are relatively prime over $K$. It
follows from Lemma \cref{\sectno.1} that the ring
$K[\tau_n]_{(f(t)/(t-\tau_n))}$ can be written as a product of fields
$K[\upsilon_i]_{h_i}$ where $K[\upsilon_i] = K[t]/(g_i)$ and $h_i(t)
=g_1(t) \cdots g_i(t)/(t-\upsilon_i)\dots g_r(t)$. It is clear that
each $g_i(t)$ is separable over $K$, and that $g_i(t)/(t-\upsilon_i)$
is separable over $K[\upsilon_i]$. In particular we have that
$K[\upsilon_i]$ is a separable field extension of $K$. It follows from
the induction hypothesis that $K[\upsilon_i]_{h_i}$ is a product of
fields that are separable over $K[\upsilon_i]$, and consequently are
separable over $K$.  Hence we can write $K_f$ as a product of
separable field extensions of $K$.

 To prove the last part of the proposition we observe that assertion
(2) of the proposition is equivalent to writing $K_f$ as a product
$K_f =\prod_{i=1}^r K_i e_i$ where $K_i$ are subfields of $K_f$ that
are separable over $K$ and $e_1, e_2, \dots , e_r$ are the primitive
idempotents of $K_f$.  Since the $e_i$ are the only primitive
idempotents in $K_f$ they are permuted by $\Sym$. The orbit of $e_i$
under this action is invariant under $\Sym$ and thus in $K$. Hence it
must be equal to $\sum_{i=1}^r e_i$, that is, $\Sym$ operates
transitively on the idempotents $e_1, e_2, \dots , e_r$. It follows
that the fields $K_i$ are isomorphic. Let $G$ be the stabilizer of one
of the idempotents $e_i$. Then $G$ operates on $K_i$ and the ring of
invariants of $K_f$ under $\Sym$ is clearly equal to the ring of
invariants of $K_i$ under $G$. The last part of the proposition thus
follows from part (1) of the proposition.
\enddemo

\proclaim{\sectno.8 Corollary}Let $L$ be the splitting field of a
separable polynomial $f(t)$ over the field $K$ and let $G$ be the
Galois group. Then $K$ is the ring of invariants under the action of
$G$ on $L$.
\endproclaim

\demo{Proof}We have a surjective map $K_f \to L$. It follows from the
proposition that $K_f =\prod_{i=1}^r K_i e_i$ where $e_1, e_2, \dots ,
e_r$ are the primitive idempotents of $K_f$, and where $L$ is
isomorphic to $K_i$ for $i=1,2, \dots , r$ for a separable field
extension $L$ of $K$. The idempotents $e_1, e_2, \dots , e_r$ map to
either $0$ or $1$ in $L$. Since $e_ie_j=0$ when $i\neq j$ we have that
exactly one $e_i$ maps to $1$ and the remaining $e_j$ to $0$. Hence
the map $K_f \to L$ factors via the projection $K_f \to K_i
e_i$. Hence we have an isomorphism $ L\to L$. It is clear that the
Galois group of $L$ is the stabilizer of $e_i$. Hence it follows from
the proposition that $K$ is the ring of invariants of $L$ under the
action of $G$.
\enddemo
 
\example{\sectno.9 Example}One should not be misled by Corollary
(2.5) and Theorem \cref{3.1} to draw too strong conclusions
about the relation between the invariants of splitting algebras and
the invariants of automorphisms of fields in the case of non-separable
extensions. To illustrate this we consider the polynomial $f(t) =t^3
-s$ over the function field $K(s)$ of the variable $s$ over a field
$K$ of characteristic $3$.  Then $K(s)_f =K(s)[\tau_1, \tau_2,
\tau_3]$ where $L= K(s)[\tau_3] =K(s) [t]/(t^3-s)$ is the splitting
field of the polynomial $f(t)$ over $K(s)$, which is purely
inseparable of degree $3$.  Let $f_2(t) =f(t)/(t-\tau_3)
=(t-\tau_3)^2$ and $K(s)_f =K(s)[\tau_3][t]/((t-\tau_2)^2)$. Hence we
have that $(\tau_i -\tau_j)^2 =0$ in $K(s)_f$ and that $\frak I
=(\tau_1 -\tau_2, \tau_1 -\tau_3,
\tau_2-\tau_3)$ is the unique maximal ideal of $K(s)_f$. We have that
the ideal $\frak I$ is stable under the action of the group $\frak
S_3$.  Hence $\frak S_3$ acts on the quotient field $K(s)_f/\frak I$
which is the splitting field of $t^3-s$ over $K(s)$.  However, we have
that $\frak S_3$ acts trivially on $K(s)_f/\frak I$ because
$\sigma\tau_i -\tau_i =\tau_{\sigma^{-1}(i)} -\tau_i$ is in $\frak I$
for all $i$.  Hence, although $K(s)$ is the ring of invariants of
$K(s)_f$ under $\frak S_3$, the group of $K(s)$-automorphisms of $L$
is trivial.  \endexample

 \example{\sectno.10 Remark} Example \cref{\sectno.9} implies that the
 extension
 $$0\to \frak I\to K(s)_f \to L\to 0$$ is non-trivial as an extension
 of $\frak S_3$-modules. More precisely, as $\frak I^2=0$, we have
 that $\frak I$ is a $1$-dimensional $L$-vector space and it is easily
 seen that $\frak S_3$ acts on it by the signum character. The
 boundary map of the above exact sequence gives an $L$-isomorphism $L
 \to H^1(\frak S_3,\frak I)$.

\endexample

\enddocument

\end